\documentclass[review]{elsarticle}

\usepackage{lineno,hyperref}
\usepackage{amssymb}
\usepackage{amsmath,amssymb,graphicx}
\usepackage{amsmath}
\usepackage{lineno,hyperref}
\usepackage{amsmath}
\usepackage{mathtools}
\usepackage{graphicx}
\usepackage{caption}
\usepackage{subcaption}
\usepackage{textcomp}
\usepackage{siunitx}
\usepackage{multirow}
\usepackage{geometry}
\usepackage{setspace}
\usepackage{color}
\usepackage{epstopdf}
\usepackage{booktabs}
\usepackage[export]{adjustbox} 
\newcommand{\rpm}{\sbox0{$1$}\sbox2{$\scriptstyle\pm$}
	\raise\dimexpr(\ht0-\ht2)/2\relax\box2 }

\modulolinenumbers[5]
\begin{document}

\begin{frontmatter}

\title{Von Neumann Stability Analysis for Multi-level Multi-step Methods}

\author[mymainaddress]{A Arun Govind Neelan\corref{mycorrespondingauthor} }
\cortext[mycorrespondingauthor]{Corresponding author}
\ead{arungovindneelan@gmail.com}
\author[mymainaddress]{\\ Department of Mechanical Engineering, Indian Institute of Technology-Madras, Chennai-600036, Tamil Nadu, India }
%

%
%
%
%

\begin{abstract}
Von Neumann stability analysis, a well-known Fourier-based method, is a widely used technique for assessing stability in numerical computations. However, as noted in "Numerical Solution of Partial Differential Equations: Finite Difference Methods" by Smith (1985, pp. 67–68), this approach faces limitations when applied to multi-level methods employing schemes with more than two levels. In this study, we aim to extend the applicability of Von Neumann stability analysis to multi-level methods. An alternative method closely related to Von Neumann stability analysis is the Approximate Dispersion Relation (ADR) analysis. In this work, we not only explore ADR analysis but also introduce various ADR analysis variants while examining their inherent limitations so that other researchers can improve the analysis before using that in their work. Furthermore, we propose an innovative strategy for reducing dissipation, optimizing it through the use of an evolutionary algorithm. Our findings demonstrate that our proposed method yields minimal errors when compared to other advection equation schemes, both in one and two spatial dimensions.

\end{abstract}

\begin{keyword}
Spectral Analysis \sep Finite difference method \sep Runge-Kutta method \sep Approximate dispersion relation analysis \sep 
\MSC[2010]
 35F61\sep 
35L50\sep 
65M06\sep 
65M08\sep 
65M12\sep 
\end{keyword}

\end{frontmatter}


\section{Introduction}
It is well known that the Lax-Richtmyer equivalence theorem \cite{https://doi.org/10.1002/cpa.3160090206} states that a consistent finite-difference method (FDM) for a linear partial differential equation for which an initial value problem is well posed is convergent if and only if it is stable. Consequently, stability stands as a pivotal criterion for evaluating the numerical solution of a consistent discretized partial differential equation. Various methods for stability analysis are employed, including Von-Neumann analysis \cite{crank1947practical}, Strong stability \cite{2007xviii, godunov1959difference}, and Eigenvalue-based stability analysis \cite{glendinning1994stability}. Among them Von-Neumann stability analysis very popular in fluid dynamics community. It was introduced briefly in a 1947 article by Crank and Nicolson \cite{crank_1996_a}, it was later formalized in a publication co-authored by John von Neumann \cite{charney_1950_numerical}. The key assumptions in Von-Neumann stability analysis involve a linear, constant-coefficient partial differential equation (PDE) with periodic boundary conditions and only two independent variables. Additionally, it prescribes that the time integration scheme should not exceed two levels, as outlined by Smith \cite{smith1985numerical}. In this work, we have overcome the last limitation and extended the Von-Neumann stability analysis to accommodate multi-level methods.

All of the aforementioned analyses can effectively predict whether the numerical error will increase or decrease. However, they do not provide advance estimates of how this error will evolve over time or iterations within the numerical solution. To assess the degree of deviation between the numerical solution and the exact solution, researchers often employ an approach known as Approximate Dispersion Analysis (ADR)~\cite{TAM1993262, LELE199216}. ADR analysis is a type of spectral analysis which quantifies the error in the numerical solution by numerical dissipation and dispersion. The ADR analysis is also similar to Von Neumann stability analysis but it also account for dispersion relation aprt from numerical amplification factor. So we can say ADR analysis is a spectral analysis build on the top of Von Neumann stability analysis. One notable advantage of ADR analysis is its capability to reduce the number of grid points required for simulations by minimizing spectral errors. Nevertheless, it's essential to recognize that not all ADR analyses found in the existing literature are universally applicable to all governing equations. In this study, we briefly outline some of the limitations associated with current ADR analysis techniques.

The earlier variant of ADR assumed time discretization is exact, we now refer to as spatial-ADR analysis. Subsequently, researchers extended this analysis to temporal discretization, assuming exact spatial derivatives, resulting in the temporal-ADR analysis. Sengupta et al. \cite{SENGUPTA20071211} recognized the significance of considering both spatial and temporal components and introduced the spatial-temporal ADR analysis, also known as global spectral analysis (GSA). While the combined spatial-temporal ADR analysis represents a more comprehensive approach for linear transport equations compared to the previous analyses. But it requires some minor refinements when applied to multi-level time integration schemes. It's important to note that GSA yields the same stability equation for a given order scheme. For instance, the stability equation for the fourth-order low-storage Runge-Kutta method is identical to that of the fourth-order strong stability preserving (SSP)-RK4 method. However, numerical simulations have shown that they exhibit distinct stability limits and slight variations in simulated results. The current analysis builds upon the foundation of the global spectral analysis by addressing some of its limitations. 

In this research, we introduce enhancements to the combined spatial-temporal ADR analysis, resulting in unique stability equations tailored to various temporal discretization methods. Notably, we've observed a precise correlation between the estimated numerical gain and the numerical gain obtained through simulations using the fast-Fourier-transform (FFT) technique. To the best of my knowledge, this is the first work present the correct analysis for multi-level, multi-step methods. 

 The paper is structured as follows:
In Section~\ref{sec:meth}, we present an extensive exploration of different ADR analysis types. Section~\ref{sec:st_adr} is dedicated to introducing our novel analysis, covering various spatial discretization and time integration schemes, including methods such as Runge-Kutta (RK) and Adam-Bashforth. We also delve into different families of RK methods, such as SSP-RK, low-storage RK, and hyperbolic Runge-Kutta. Test cases are detailed in Section~\ref{sec:result}, offering practical insights. In Section~\ref{sec:RK6L4R2}, we present an optimized scheme which minimize dissipation error. We also address the limitations of our analysis in Section~\ref{sec:con}. To facilitate reader comprehension, we've thoughtfully included MATLAB scripts for certain discretization methods, accessible via the link provided in the report.  

\section{Review of ADR analysis}
\label{sec:meth}
The ADR analysis has undergone a significant evolution over several years, thanks to the contributions of numerous researchers~\cite{SENGUPTA20071211,https://doi.org/10.1111/j.2153-3490.1972.tb01547.x,doi:10.1146/annurev.fl.06.010174.001433}. To distinguish between these contributions, we categorize them into three distinct types. The first approach is primarily geared toward spatial discretization, the second one pertains to temporal discretization, and the third type addresses the intricacies of combined spatial and temporal discretization.
\subsection{Spatial ADR analysis}
Spatial ADR analysis~\cite{https://doi.org/10.1111/j.2153-3490.1972.tb01547.x,doi:10.1146/annurev.fl.06.010174.001433} is a popular ADR analysis extensively used in ~\cite{tan2022two, DESHPANDE2021110157}.
For spatial ADR analysis, we write spatial discretization terms in the spectral plane. The general form of the first derivative term for the finite difference scheme is written as
\begin{equation}\label{eq:gfd}
	\left.\frac{d u}{d x}\right|_{i} = \frac{1}{h} \sum_{j=-N}^{M}c_j u_{i+j}.
\end{equation} 
where, $h$ is the grid size.
The Fourier representation of the signal is
\begin{equation}\label{key}
	u(x) = \int \tilde{U} (k) \exp(ikx) \ dk.
\end{equation}
Here, $i =\sqrt{-1}$,  $\tilde{U}$ is amplitude of frequency. Writing the \eqref{eq:gfd} in spectral plane gives,
\begin{equation}\label{key}
	ik\tilde{U} (k) =\frac{1}{\Delta x} \left(\sum_{j=-N}^{M}c_j \exp(ik j \Delta x)\right)\tilde{U}(k).
\end{equation}
The modified wave number is
\begin{equation}\label{eq:k_num}
	k_{num} =\frac{-i}{\Delta x} \left(\sum_{j=-N}^{M}c_j \exp(ik j \Delta x)\right).
\end{equation}
To minimize spectral errors, it is essential to ensure that the modified wave number in equation \eqref{eq:k_num} aligns with the physical wave number. The real component of the modified wave number characterizes the numerical dissipation introduced by the discretization, while the imaginary component accounts for the numerical dispersion. The first-order backward in space discretization is
\begin{equation}\label{key1}
	\left.\frac{d u}{d x}\right|_{i} = \frac{u_{i}-u_{i-1}}{\Delta x} + \mathcal{O}(\Delta x).
\end{equation} 
The modified wave number equation of the \eqref{key1} is
\begin{equation}\label{key}
	\frac{k}{k_{eq}} = \frac{\sin(kh)}{kh}- i\frac{1-\cos(kh)}{kh}.
\end{equation}
\par The initial segment of the equation pertains to numerical dispersion, while the subsequent segment relates to numerical dissipation. 
Spatial ADR analysis proves effective in predicting numerical dissipation, particularly in scenarios where time integration remains exact. However, its limitation arises from the fact that it solely addresses spatial aspects and doesn't account for temporal terms. Consequently, it may not be well-suited for solving transport equations that encompass both spatial and temporal components. In the context of transport equations, the choice of time steps becomes pivotal in stability analysis, as evidenced by the CFL number ($c\frac{\Delta t}{h}$).
Spatial ADR analysis, unfortunately, overlooks the influence of $\Delta t$, and consequently, it does not consider the impact of the CFL number. Nevertheless, it's worth noting that this analysis has played a significant role in the evolution of ADR analysis methodologies.
\subsection{Temporal ADR analysis} 
Similar to spatial discretization, ADR for temporal discretization is presented here. More details about the temporal ADR analysis and optimzing scheme based on temporal ADR can be found in~~\cite{BOGEY2004194, HU1996177}.
\begin{eqnarray}\label{key}
	\begin{aligned}
		& u^{n+1}= u^n +\sum_{i=1}^{p}w_ik_i,
		& k_i=\Delta t  F\bigg(U^n+\sum_{j=1}^{i-1}\beta_{i,j}k_j\bigg),\qquad i = 1, 2,..., p.
	\end{aligned}
\end{eqnarray}
Here $w_i$ and $\beta_{i,j}$ are the weights of the multi-step method. Writing the multi-step method in the spectral plane, we have
\begin{equation}\label{key}
	\tilde{U}^{n+1}= \tilde{U}^{n}\left(1+\sum_{j=1}^{p}w_j(ick\Delta t)^j\right),
\end{equation}
where, $\tilde{U}$ represents the frequency-domain representation of the signal and $c$ is speed of linear convection equation.  The numerical amplification factor of linear convection equation is~\cite{BOGEY2004194, HU1996177}
\begin{equation}\label{key}
	G_{theory}=\frac{\tilde{U}^{n+1}}{\tilde{U}^{n}}=1+\sum_{j=1}^{p}\left(w_j(ick\Delta t)^j\right).
\end{equation}

\par
\begin{figure}[h]
	\begin{subfigure}[b]{0.5\textwidth}
		\includegraphics[width=\linewidth]{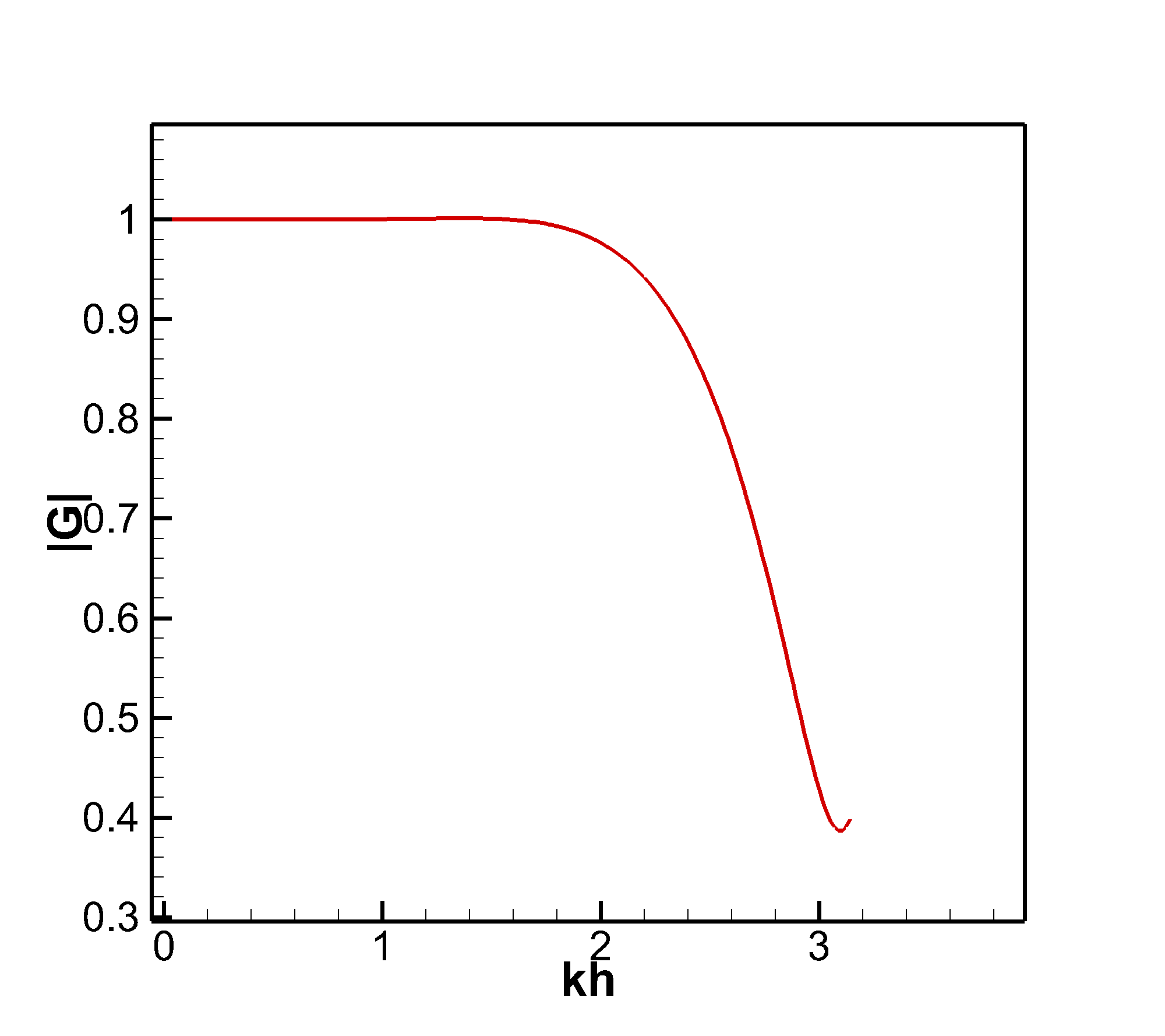}
		\begin{center}
			\caption{$|G|=1$ plot for RK45 }
			\label{fig:t_G_RK451}
		\end{center}
		
	\end{subfigure}%
	\begin{subfigure}[b]{0.5\textwidth}
		\includegraphics[width=\linewidth]{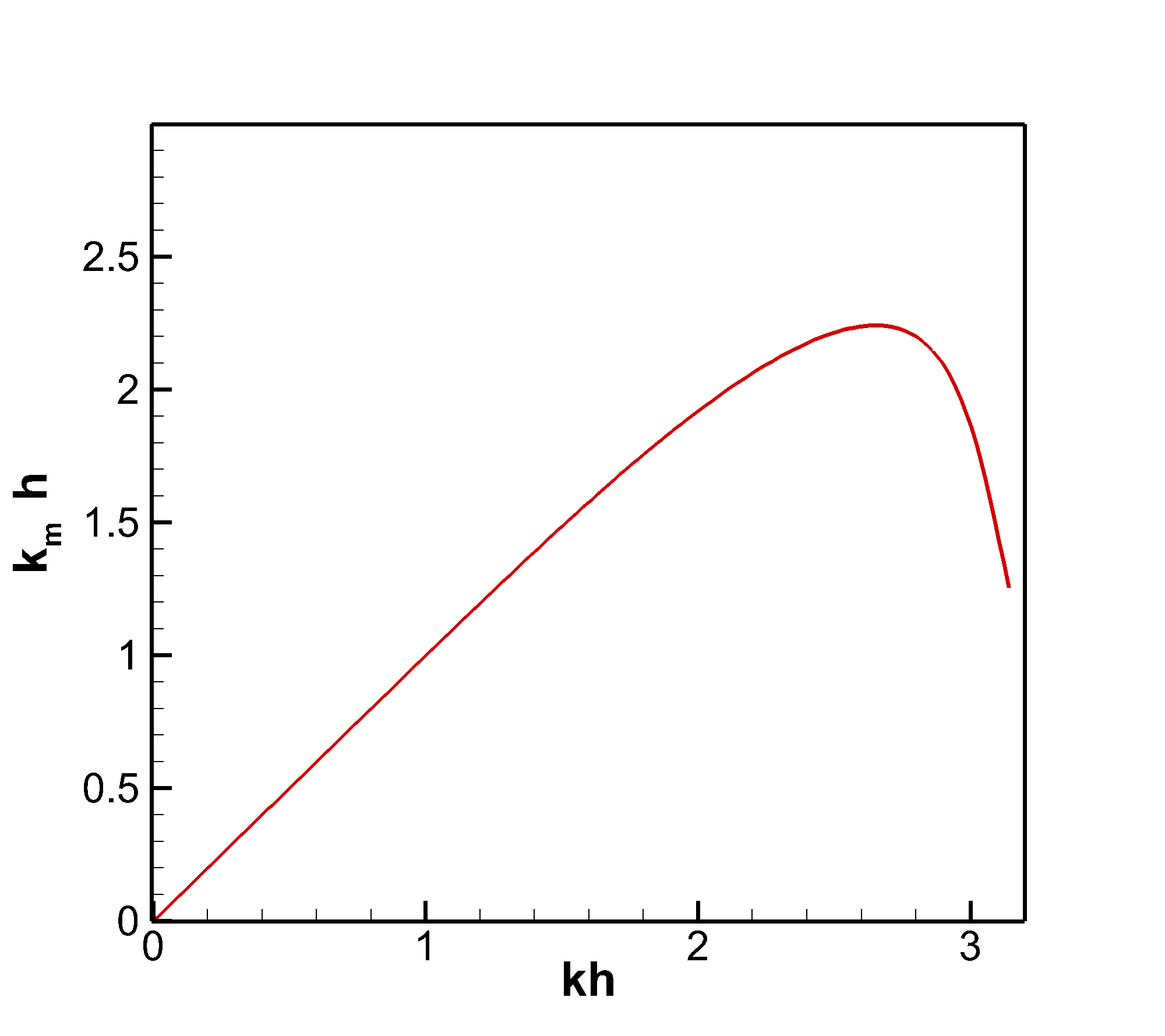}
		\caption{Modified wave number plot for RK45  }
		\label{fig:t_w_RK45}
	\end{subfigure}
	\begin{center}
		\caption{Temporal ADR plot of RK45}
	\end{center}
	\label{fig:tRK45}
\end{figure}

 In order to demonstrate the constraints associated with temporal ADR analysis, a novel fourth-order, five-stage Runge-Kutta (RK45) method was developed. This method was designed to uphold a specific constraint on the numerical amplification factor ($|G| \le 1$) for all values of $kh$.
The RK45 method we developed to illustrate the limitation is
\begin{eqnarray}
	\begin{aligned}
		& u^{(1)}= u^{(n)}+0.18856248 \Delta t  F(u^{(n)}), \\
		& u^{(2)}= u^{(n)}+\frac{1}{4} \Delta t  F (u^{(1)}),\\
		&u^{(3)}= u^{(n)}+\frac{1}{3}\Delta t  F (u^{(2)}),\\
		& u^{(4)}= u^{(n)}+\frac{1}{2} \Delta t  F (u^{(1)}),\\
		&u^{(5)}= u^{(n)}+\frac{1}{1}\Delta t  F (u^{(1)}),\\	
	\end{aligned}
\end{eqnarray}
The theoretical amplification factor is
$$G^{T}_{theory} = 1+\frac{ (i \omega \Delta t)}{1}+ \frac{ (i \omega \Delta t)^2}{2!}+\frac{ (i \omega \Delta t)^3}{3!}+\frac{ (i \omega \Delta t)^4}{4!}+0.18856248\times\frac{1}{24}{ (i \omega \Delta t)^5}.$$
The $|G_{theory}|$ plot for the RK45 scheme is shown in figure~\ref{fig:t_G_RK451} and the modified wave number plot is displayed in figure~\ref{fig:t_w_RK45}. 
\begin{figure}[h!]
	\centering
	\includegraphics[width=0.5\linewidth]{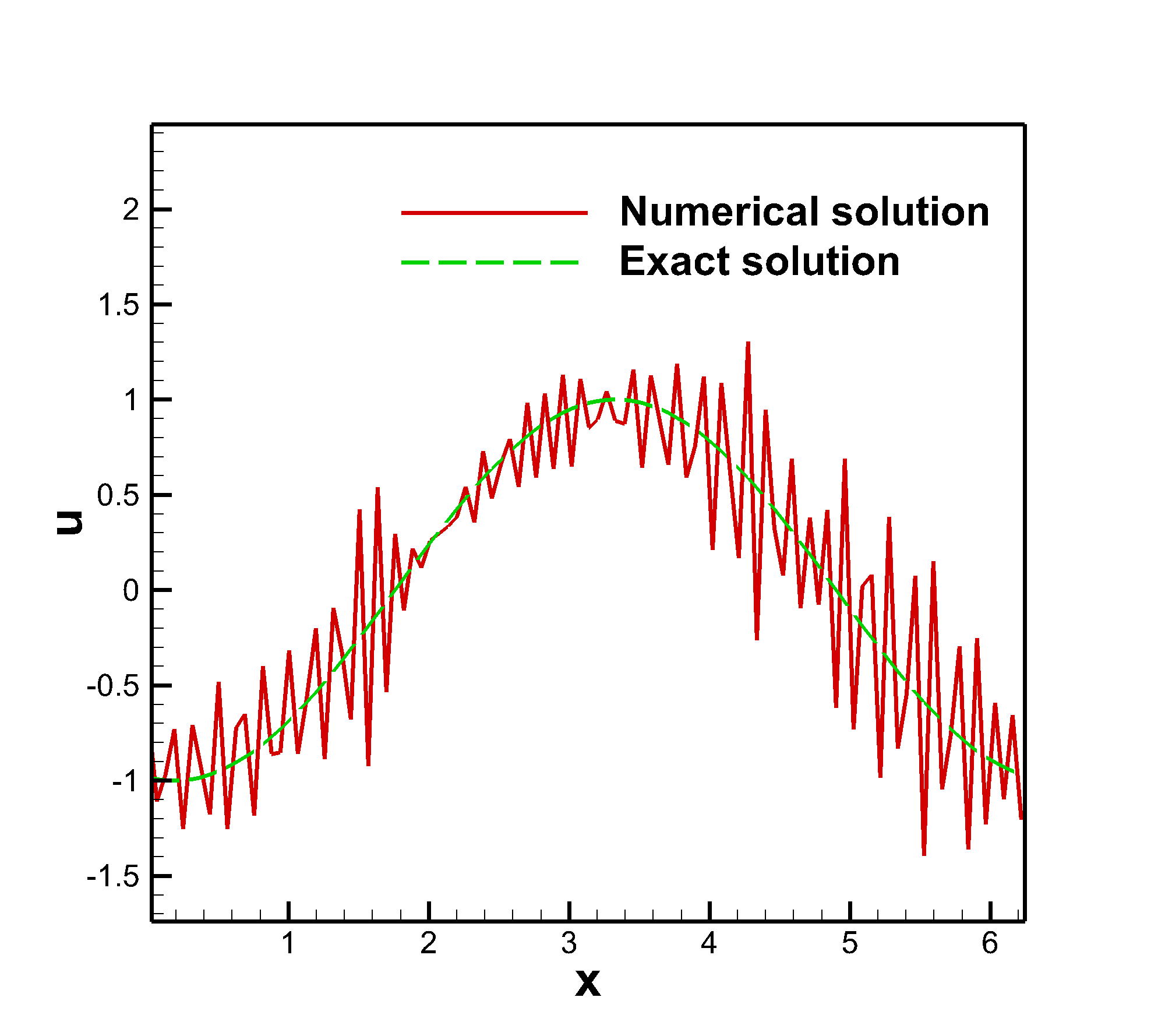}
	\caption{Solution of Linear convection equation with CFL = 4 at 1.7~s}
	\label{fig:rk45wave}
\end{figure}
From this, it is clear that the magnitude of $G_{theory}$ is below one for all wave numbers so this scheme is supposed to be stable regardless of the solution setup. We shall test this using the  one-dimensional linear convection equation. The linear convection equation is 
\begin{equation}\label{key}
	\frac{\partial u}{\partial t}+c\frac{\partial u}{\partial x} =0,
\end{equation} 
In this scenario, the variable $c$ represents the signal's speed, and it is determined using the RK45 method for temporal integration. The spatial discretization is executed using a backward-in-space scheme across 101 grid points spanning the domain [0, $2\pi$]. A CFL number of 4 is employed, and the simulation runs up to a flow time of 1.7 seconds. The resulting solution of the equation is presented in Figure~\ref{fig:rk45wave}. Notably, oscillations are evident in the solution, and these oscillations tend to amplify as time progresses. 

This observation contradicts the expectations based on temporal ADR analysis, which suggests stability over time steps. However, the presence of oscillations in the results can be attributed to the fact that the current analysis is limited to temporal discretization. It is more suitable for ordinary differential equations that solely involve a time-dependent term. To address these limitations, a combined analysis, as described in the following section, is proposed to overcome these challenges.	

\subsection{Space-time ADR analysis}

The foundation of this analysis was initially laid out in the publication by Sengupta et al. ~\cite{SENGUPTA20071211}. 
However, it's essential to acknowledge that this method does have certain constraints when applied to multi-level methods, and these limitations will be discussed shortly. The current approach draws inspiration from and builds upon the foundation laid by the method mentioned above. 

We revisit the linear convection equation 
\begin{equation}\label{key}
	\frac{\partial u}{\partial t}+c\frac{\partial u}{\partial x}=0.
\end{equation}
The first step in the space-time spectral analysis is writing the equation in the spectral plane,
\begin{equation}\label{key}
	u(x,t)=\int\int \tilde{U}(k,\omega) \exp(i(kx-\omega t)) dk d\omega.
\end{equation}
The dispersion relation,
$
\omega = ck,
$ is obtained 
by substituting standing wave solution, $u(x,t)= \tilde{U}(k,\omega) \exp(i(kx-\omega t))$ in linear convection equation. This determines the signal propagation speed of different wave numbers. When several wave numbers are present in the signal, different wave numbers may propagate at different velocities so group velocity is derived. The group velocity is defined as
\begin{equation}\label{key}
	V_{g,theory}=\frac{\partial \omega_{theory}}{\partial k}.
\end{equation}
The growth or decay of the signal over time is calculated by gain. The numerical gain is
\begin{equation}\label{eq:G1}
	G_{theory} = \frac{\tilde{U}(k,t+\Delta t)}{\tilde{U}(k,t)} = \exp (-ikh\, \text{CFL}).
\end{equation}
where $CFL =c\frac{\Delta t}{h}$, $h$ is grid size, $c$ is speed of the signal and $\Delta t$ is time-step. The numerical phase shift from theory is 
\begin{equation}\label{key}
	\tan (\beta)  = \left[\frac{(G_{theory})_{Imag}}{(G_{theory})_{Real}}\right],
\end{equation}
and the numerical phase speed is
\begin{equation}\label{key}
	\frac{c_{theory}}{c}=\frac{\beta}{kc\Delta t}=\frac{1}{kh\ \text{CFL}} \tan^{-1}\left[\frac{(G_{theory})_{Imag}}{(G_{theory})_{Real}}\right],
\end{equation}
The  group velocity can be calculated using
\begin{equation}\label{key}
	V_{g,theory}=\frac{\partial \omega_{theory}}{\partial k}.
\end{equation}
Using $\omega_{theory}=kc$, further simplification gives
\begin{equation}\label{eq:Vg}
	\frac{V_{g,theory}}{c}=\frac{1}{\text{CFL}}\frac{d \beta}{d (kh)}.
\end{equation}

The numerical amplification factor of the  $p^{th}$ order RK method given in \eqref{eq:rkg} is\cite{SENGUPTA2020109310}
\begin{equation}\label{eq:sen}
	G_{thoery} = 1+\sum_{j=1}^p (-1)^ja_jA_i^j
\end{equation} 
Here, we have $a_j = 1/j!$, and $A_i$ represents the Fourier-transform of the spatial discretization method employed. It's important to note that the numerical amplification factor plot described in \eqref{eq:sen} remains consistent for Runge-Kutta methods of a specific order. However, within a given order, there exists a broad range of RK methods, each with distinct stability characteristics. For instance, classical RK4 and SSPRK4 exhibit different simulation outcomes and stability limits, despite both methods yielding the same stability equation through analysis.  

\section{Von Neumann stability analysis}
\label{sec:st_adr}
The existing spatial-temporal ADR analysis gives same stability equation for a given order of time integration scheme used. This can be easily modified by incorporating multi-level spectral analysis to the scheme is presented in this section. In this work, we restrict our scope only numerical gain and we will not study about dispersion relation. So we call in analysis Von Neumann stability analysis not ADR analysis. Extending this work to ADR analysis is straight forward.   

\subsection{Forward in time and backward in space (FTBS) scheme }  
\label{sec:ftbs}
FTBS scheme for linear convection equation is
\begin{equation}\label{eq:FTBS}
	u^{n+1}=u^n-\frac{c\Delta t}{h}(u^n_{i}-u^n_{i-1}).
\end{equation}
Converting the above equation in the spectral plane, we get
\begin{equation}\label{key}
	\tilde{U}^{n+1}=\tilde{U}^{n} - CFL \left(\tilde{U}^{n}-\tilde{U}^{n} \exp(-ikh)\right),
\end{equation}
upon simplification. we yield
\begin{equation}\label{eq:g_ftbs}
	G_{FTBS} = 1 -CFL (1-\exp (-ikh)),
\end{equation}
Substituting CFL = 1 into \eqref{eq:g_ftbs} yields $G_{\text{FTBS}} = 1$. When $|G|<1$, the initial condition experiences damping, and conversely, when $|G| > 1$, it exhibits amplification. At $|G| = 1$, there is no dissipation error. In the absence of errors in the numerical solution, the numerical result matches the analytical solution. This is why FTBS can provide the exact solution for the linear convection equation when solved with CFL = 1.

	\begin{figure}[h]
		\centering
		\includegraphics[width=.5\linewidth]{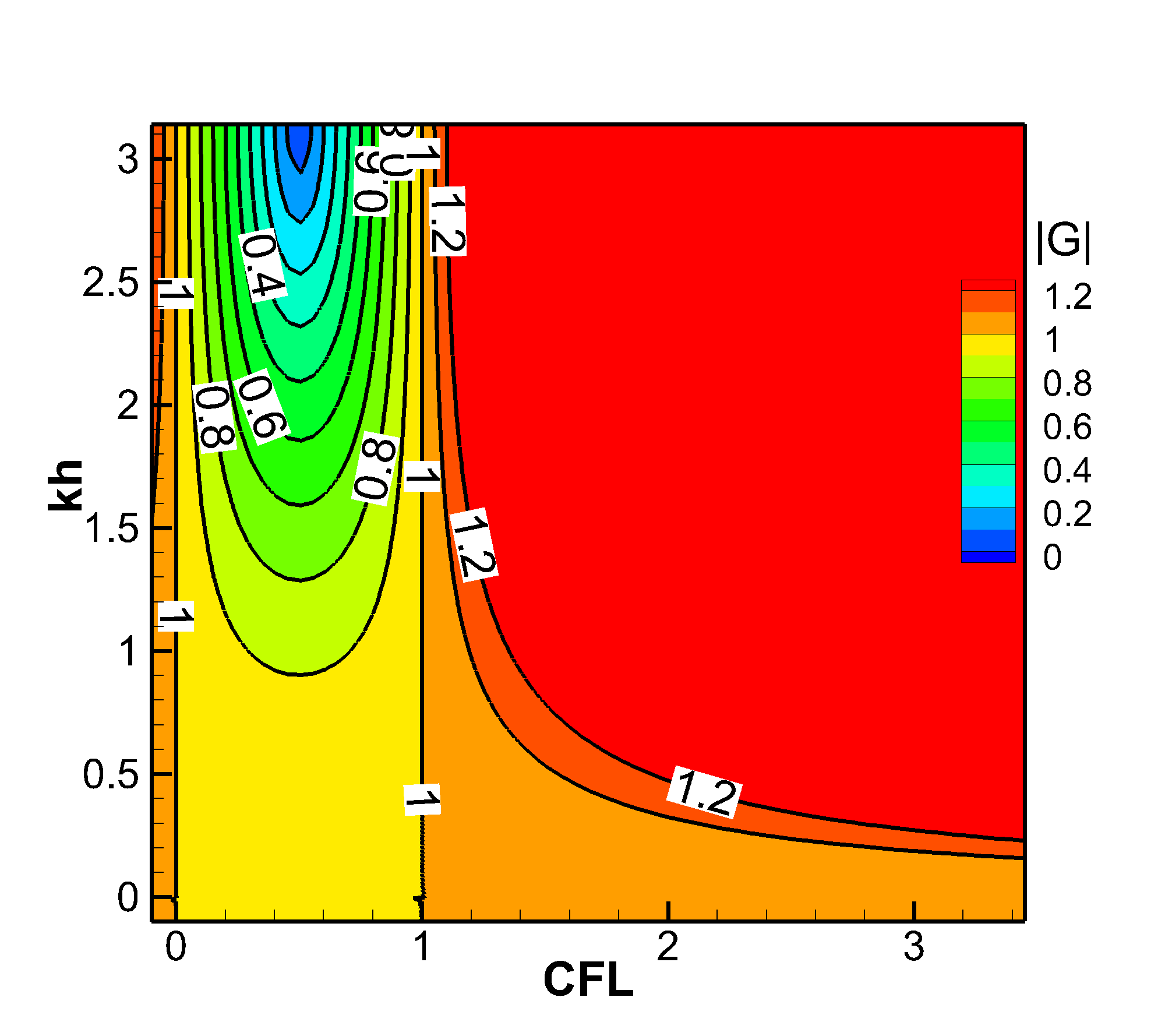}
		\caption{$|G|$ plot for FTBS }
		\label{fig:G_FTBS_1}
	\end{figure}%

Figure~\ref{fig:G_FTBS_1} displays the gain plot magnitude of the FTBS scheme. Notably, when $G = 1$, it forms a straight line perpendicular to the x-axis. This signifies that the discretization is capable of resolving any wave number when CFL = 1. It's important to emphasize that the FTBS discretization is employed to solve the linear convection equation across the domain [0, $2\pi$], utilizing 101 grid points and an initial condition of $u(x,0)=\sin(x)$. Furthermore, it's noteworthy that the FTBS scheme with CFL = 0.2 and 0.8 yields theoretical gains of 0.9996 - 0.0126i and 0.9984 - 0.0502i, respectively, which closely align with the gains determined through numerical calculations. These numerical gains are computed by taking the Fast Fourier Transform (FFT) at the current time and dividing it by the FFT at the previous time in the simulation. The readers are encouraged to use the Matlab code given in the \href{https://github.com/AGN000/Von_Neumann_Stability_matlab_code}{link} to study the stability behavior of simulated result and theoretical prediction 

\subsection{SSP-RK3-L2R1 scheme}
The SSP-RK3 method is a popular RK method designed for better stability~\cite{doi:10.1137/S003614450036757X}. Here, we have discretized the linear convection equation using the three-stage, third-order SSPRK3 method for temporal discretization and the third-order four-point ($i-2$, $i-1$, $i$ and $i+1$) finite difference scheme for space. The formulation used is
\begin{align}
	& u^{(1)} = u^n - \frac{c\Delta t}{6\Delta x}\left( u^{n}_{i-2}-6u^{n}_{i-1}+3u^{n}_{i}+2u^{n}_{i+1}\right), \\ \nonumber
	& u^{(2)} = \frac{3}{4}u^n +\frac{1}{4}u^{(1)}- \frac{c\Delta t}{4\times 6\Delta x}\left( u^{(1)}_{i-2}-6u^{(1)}_{i-1}+3u^{(1)}_{i}+2u^{(1)}_{i+1}\right),  \\ \nonumber
	& u^{n+1}=\frac{1}{3}u^n+\frac{2}{3}u^{(2)} - \frac{2c\Delta t}{3\times6\Delta x}\left( u^{(2)}_{i-2}-6u^{(2)}_{i-1}+3u^{(2)}_{i}+2u^{(2)}_{i+1}\right).
\end{align}
Applying the Fourier-Laplace transform, we obtain
\begin{align}
	&\tilde{U}^{(1)} = \tilde{U} - \frac{CFL}{6}\left(\tilde{U}\exp(-2ikh)-6\tilde{U}\exp(-ikh)+ 3\tilde{U}+2\tilde{U}\exp(ikh)\right), \\ \nonumber
	& \tilde{U}^{(2)} = \frac{3}{4}\tilde{U} +\frac{1}{4}\tilde{U}^{(1)}- \frac{CFL}{4\times6}\left(\tilde{U}^{(1)}\exp(-2ikh)-6\tilde{U}^{(1)}\exp(-ikh)+ 3\tilde{U}^{(1)}+2\tilde{U}^{(1)}\exp(ikh)\right) , \\ \nonumber
	&\tilde{U}^{n+1}=\frac{1}{3}\tilde{U}+\frac{2}{3}\tilde{U}^{(2)} - \frac{2CFL}{3\times6}\left(\tilde{U}^{(2)}\exp(-2ikh)-6\tilde{U}^{(2)}\exp(-ikh)+ 3\tilde{U}^{(2)}+2\tilde{U}^{(2)}\exp(ikh)\right).
\end{align}
Upon simplification, we have
\begin{align} \label{eq:G_ssprk3-BS}
	&	G^{(1)} = 1 - \frac{CFL}{6}\left(\exp(-2ikh)-6\exp(-ikh)+ 3+2\exp(ikh)\right), \\ \nonumber
	&G^{(2)} = \frac{3}{4} +\frac{1}{4}G^{(1)}- \frac{CFL}{4\times6}\left(G^{(1)}\exp(-2ikh)-6G^{(1)}\exp(-ikh)+ 3G^{(1)}+2G^{(1)}\exp(ikh)\right), \\ \nonumber
	&G_{SSPRK3-L2R1}=\frac{1}{3}+\frac{2}{3}G^{(2)} - \frac{2CFL}{3\times6}\left(G^{(2)}\exp(-2ikh)-6G^{(2)}\exp(-ikh)+ 3G^{(2)}+2G^{(2)}\exp(ikh)\right).
\end{align}
Please note that \eqref{eq:G_ssprk3-BS} provides the gain characteristics of the SSPRK3-L2R1 scheme.

\begin{figure}
	\centering
	\includegraphics[width=.5\linewidth]{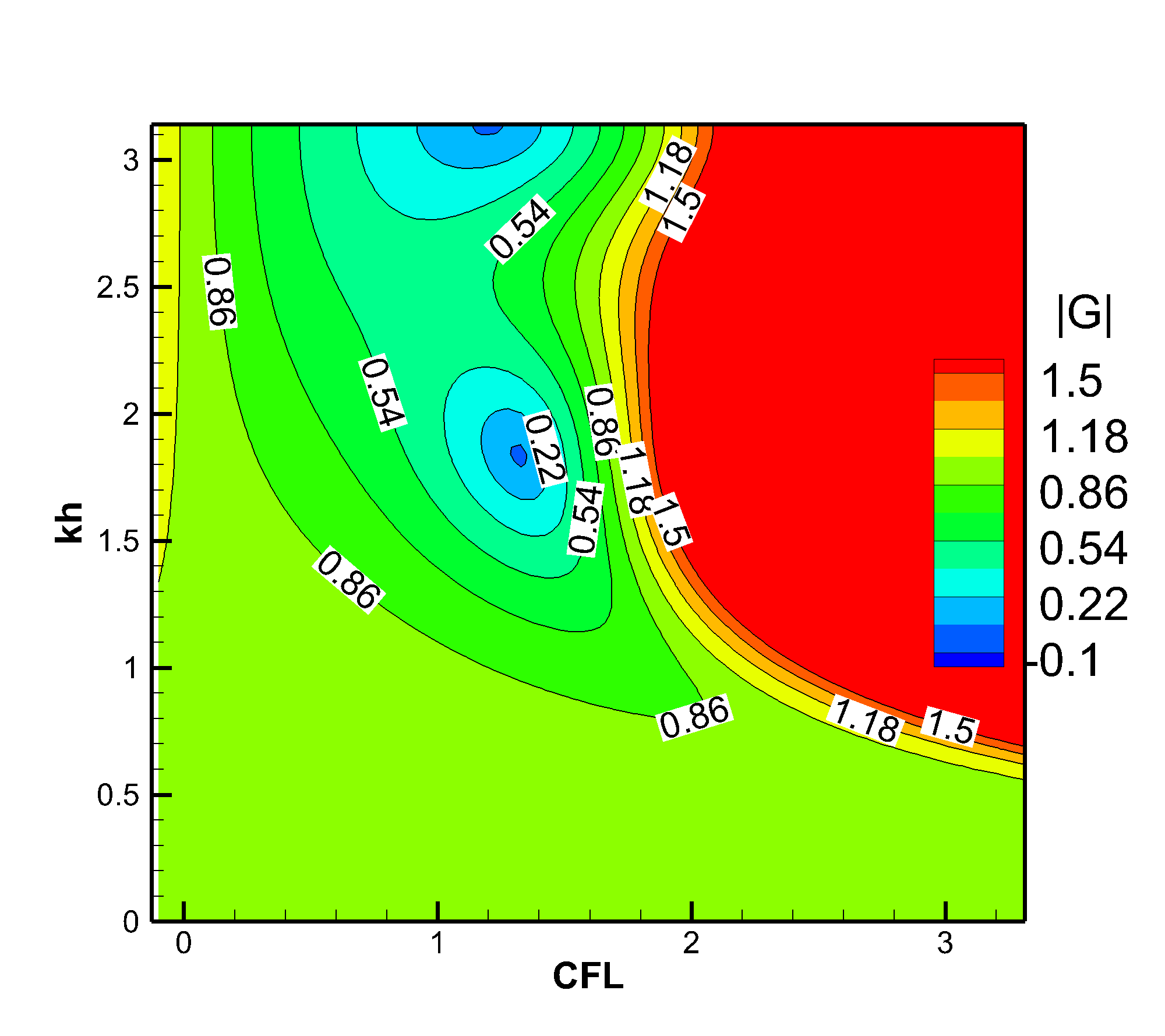}
	\caption{$|G|$ plot for SSPRK3-L2R1 }
	\label{fig:G_SSPRK3_L2R1}
\end{figure}

Figure~\ref{fig:G_SSPRK3_L2R1} illustrates the theoretical gain profiles for the SSPRK3-L2R1 discretization applied to the linear convection equation. When subject to the test condition outlined in section~\ref{sec:ftbs}, using CFL = 0.5, the numerical gain and theoretical gain for this discretization are calculated as 0.9985 - i 0.0314. This scheme demonstrates robustness at higher CFL numbers with minimal stability concerns. It maintains stability for CFL values up to 1.6 across all wave numbers. Similar to other schemes, the numerical dispersion decreases when both $kh$ and the CFL number are reduced during simulation. The SSPRK3-L2R1 scheme consistently delivers accurate results when employed with low values of $kh$ and a CFL value of $\le$ 1.5. 

\subsection{Adams-Bashforth method}
\label{subsec:AB}
Adams-Bashforth method is a linear multi-step method used to solve differential equations and transport equations. The discretizing linear convection equation using the Adams-Bashforth method for time and second-order backward scheme in space gives
\begin{equation}\label{key}
	u^{n+2*}_i=u^{n+1}_i-\frac{3CFL}{2\times2}\left(u^{n+1}_{i-2}-4u^{n+1}_{i-1}+3u^{n+1}_{i}\right)+\frac{CFL}{2\times2}\left(u^{n}_{i-2}-4u^{n}_{i-1}+3u^{n}_{i}\right).
\end{equation}
The corrector step is
\begin{equation}\label{key}
	u^{n+2}_i=u^{n+1}_i-\frac{CFL}{2}\left(u^{n+2*}_{i-2}-4u^{n+2*}_{i-1}+3u^{n+2*}_{i}\right)-\frac{CFL}{2}\left(u^{n+1}_{i-2}-4u^{n+1}_{i-1}+3u^{n+1}_{i}\right).
\end{equation}
Applying the Fourier-Laplace transform gives,
\begin{align}\label{key}
	&U^{n+2*}=U^{n+1}-\frac{3CFL}{2\times2}\left(3U^{n+1}-4U^{n+1}\exp(-ikh)+U^{n+1}\exp(-2ikh)\right)\\ \nonumber
	&+\frac{CFL}{2\times2}\left(3U^{n}-4U^{n}\exp(-ikh)+U^{n}\exp(-2ikh)\right),\\
	&U^{n+2}=U^{n+1}-\frac{CFL}{2\times2}\left(3U^{n+2*}-4U^{n+2*}\exp(-ikh)+U^{n+2*}\exp(-2ikh)\right)\\\nonumber
	&-\frac{CFL}{2\times2}\left(3U^{n+1}-4U^{n+1}\exp(-ikh)+U^{n+1}\exp(-2ikh)\right).
\end{align}
Writing in terms of numerical amplification factor gives
\begin{equation}\label{eq:2G_AB}
	G^{(1)}=1-\frac{3CFL}{2\times2}\left(3-4\exp(-ikh)+\exp(-2ikh)\right)+\frac{CFL}{2G^{(1)}}\left(3-4\exp(-ikh)+\exp(-2ikh)\right),
\end{equation}
Where $G^{(1)}=\frac{U^{n+2*}}{U^{n+1}} = \frac{U^{n+1}}{U^{n}}.$  It has two numerical gain values. Solving for $G^{(1)}$ and substitute that in the corrector equation gives
\begin{equation}\label{key}
	G=1-\frac{G^{(1)}.CFL}{2\times2}\left(3-4\exp(-ikh)+\exp(-2ikh)\right)-\frac{CFL}{2\times2}\left(3-4\exp(-ikh)+\exp(-2ikh)\right).
\end{equation}

\begin{figure}[ht!]
	\begin{subfigure}[b]{0.5\textwidth}
		\includegraphics[width=\linewidth]{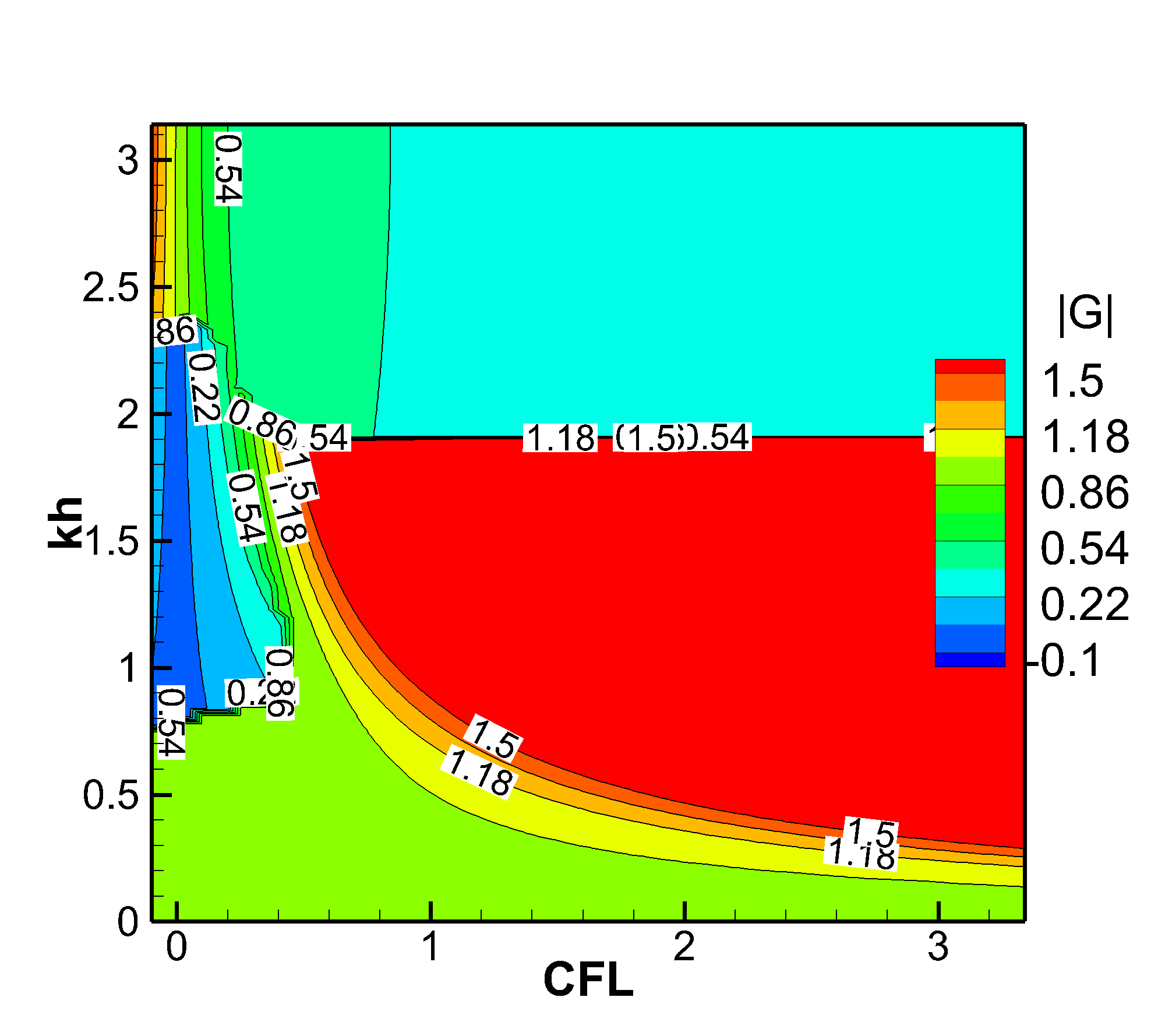}
		\caption{$|G_1|$ plot }
		\label{fig:G_AB_BS}
	\end{subfigure}%
	\begin{subfigure}[b]{0.5\textwidth}
		\includegraphics[width=\linewidth]{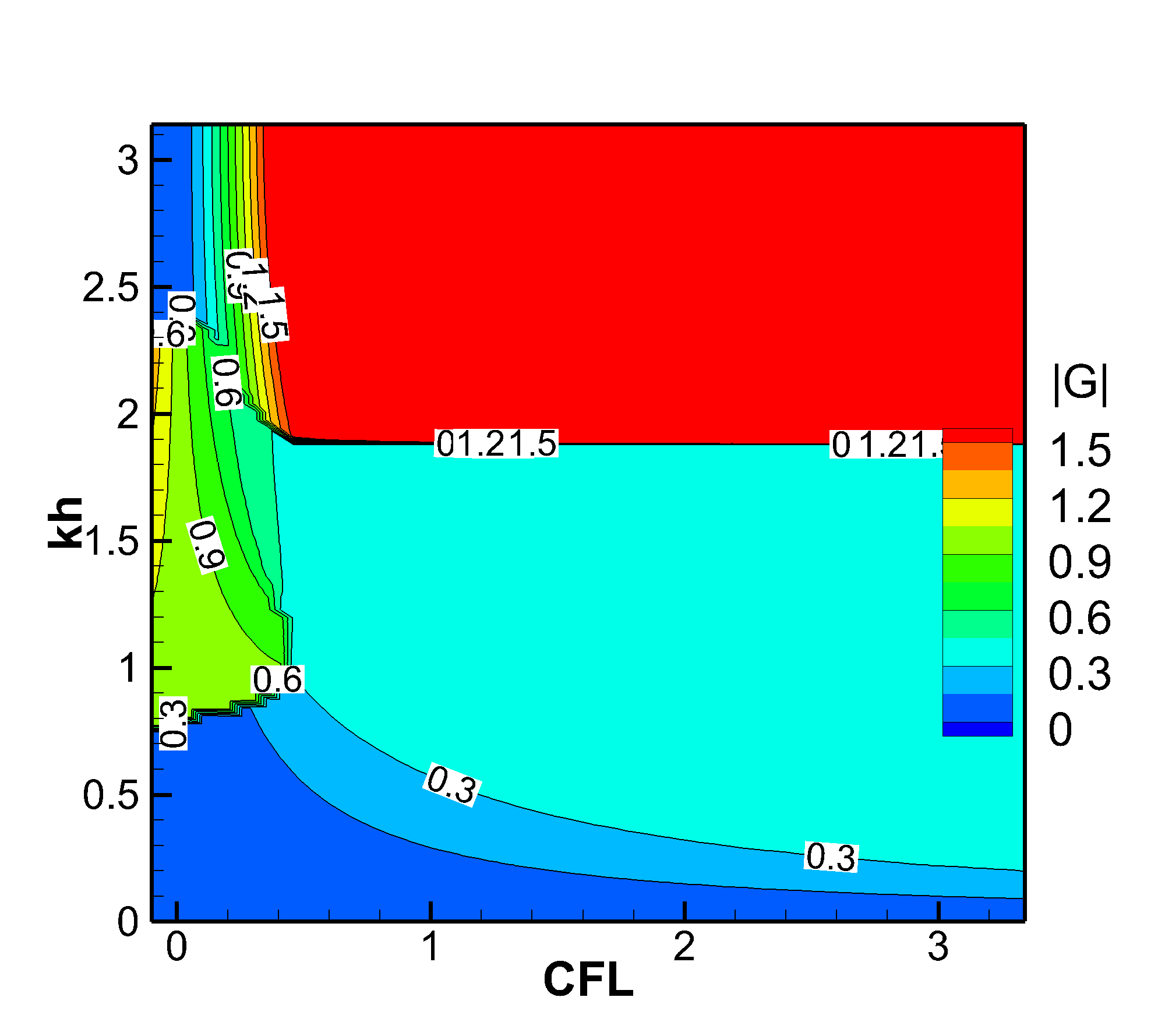}
		\caption{$|G_2|$ plot }
		\label{fig:G2_AB_BS2}
	\end{subfigure}
	\caption{$G$ plot of Adams-Bashforth method-BS2 scheme }
	\label{fig:GAB_BS}
\end{figure}

The numerical analysis yields two values for $G$, specifically denoted as $G_1$ and $G_2$. To solve the convection equation, this discretization is applied with an initial condition of $u(x,0) = \sin(x)$ and periodic boundary conditions spanning the domain [0, 2$\pi$]. The domain is divided into 100 cells of equal size, and the solution is computed up to a flow time of 1 second using a CFL number of 0.4. The numerical gain determined from the simulation results is calculated as 0.9989 - i 0.0251. In contrast, the theoretical analysis provides two values for the gain: 0.9997 - i 0.0252 and 0.0003 - i 0.0126. The first value appears reasonable, while the second one is considered unrealistic~~\cite{SENGUPTA2017182,SENGUPTA201741}.
It's noteworthy that the numerical gain obtained from the simulation aligns perfectly with one of the values obtained from the theoretical estimate. Visual representations of $G_1$ and $G_2$ for this discretization are presented in Figure~\ref{fig:G_AB_BS} and Figure~\ref{fig:G2_AB_BS2}, respectively.

While the Adam-Bashforth method can attain a specific order of accuracy with fewer function evaluations compared to Runge-Kutta (RK) methods, its popularity is relatively limited, possibly due to the following reasons:

\begin{enumerate}
	\item When dealing with transport equations, we typically compute the time step based on a given CFL number, resulting in a variable time step that adapts to flow properties. Achieving higher-order accuracy with the Adam-Bashforth method requires its adaptation to non-uniform time steps, which can be computationally expensive. Therefore, the Adam-Bashforth method might serve as a viable alternative to the RK method when a constant time step is employed to solve a transport problem.
	
	\item As illustrated in Figure~\ref{fig:G_AB_BS}, it becomes evident that the stability region of the Adam-Bashforth method is considerably smaller than that of the RK method. This limitation not only discourages the use of the Adam-Bashforth method for linear convection equations but also for transport equations in general.
	
\end{enumerate}

\subsection{RK6L4R2 scheme}
\label{sec:RK6L4R2}
In the FTBS (Forward-Time Backward-Space) scheme, we observed that when CFL (Courant-Friedrichs-Lewy) number equals 1, we reach a condition where $|G|=1$, indicating an undamped solution. This occurs because at CFL = 1, the stability region forms a straight line perpendicular to the x-axis, encompassing all wavenumbers. Consequently, the solution remains free of dissipation effects across all wavenumbers. 

In our pursuit of finding a less-dissipative scheme, we devised an evolutionary algorithm aimed at identifying a scheme where the condition $G = 1$ forms a nearly straight line over a range of $kh$ values and remains perpendicular to the x-axis. Through this algorithm, we discovered the RK6L4R2 scheme. This scheme leverages the sixth-order six-stage low-storage Runge-Kutta method for time integration and employs a seven-point ($i-4$ to $i+2$) sixth-order first derivative term for spatial discretization. The discretized form of this scheme is as follows:
\begin{eqnarray}
	\begin{aligned} \label{eq:RK6l4r2}
		& u^{(1)} = u^0 -\frac{CFL}{360}\left({u}^0_{i-4}-8{u}^0_{i-3}+30{u}^0_{i-2}-80{u}^0_{i-1}+35{u}^0_{i}+24{u}^0_{i+1}-2{u}^0_{i+2}\right), \\ \nonumber
		& u^{(2)} = u^0 -\frac{CFL}{300}\left({u}^{(1)}_{i-4}-8{u}^{(1)}_{i-3}+30{u}^{(1)}_{i-2}-80{u}^{(1)}_{i-1}+35{u}^{(1)}_{i}+24{u}^{(1)}_{i+1}-2{u}^{(1)}_{i+2}\right), \\ \nonumber
		& u^{(3)} = u^0 -\frac{CFL}{240}\left({u}^{(2)}_{i-4}-8{u}^{(2)}_{i-3}+30{u}^{(2)}_{i-2}-80{u}^{(2)}_{i-1}+35{u}^{(2)}_{i}+24{u}^{(2)}_{i+1}-2{u}^{(2)}_{i+2}\right), \\ \nonumber 
		& u^{(4)} = u^0 -\frac{CFL}{180}\left({u}^{(3)}_{i-4}-8{u}^{(3)}_{i-3}+30{u}^{(3)}_{i-2}-80{u}^{(3)}_{i-1}+35{u}^{(3)}_{i}+24{u}^{(3)}_{i+1}-2{u}^{(3)}_{i+2}\right), \\ \nonumber
		& u^{(5)} = u^0 -\frac{CFL}{120}\left({u}^{(4)}_{i-4}-8{u}^{(4)}_{i-3}+30{u}^{(4)}_{i-2}-80{u}^{(4)}_{i-1}+35{u}^{(4)}_{i}+24{u}^{(4)}_{i+1}-2{u}^{(4)}_{i+2}\right), \\ \nonumber
		& u^{n+1} = u^0 -\frac{CFL}{60}\left({u}^{(5)}_{i-4}-8{u}^{(5)}_{i-3}+30{u}^{(5)}_{i-2}-80{u}^{(5)}_{i-1}+35{u}^{(5)}_{i}+24{u}^{(5)}_{i+1}-2{u}^{(5)}_{i+2}\right) ,\nonumber
	\end{aligned}
\end{eqnarray}
where $u^0=u^n.$ Applying  the Fourier-Laplace transform  gives
\begin{eqnarray}
	\begin{aligned} 
		U^{(1)} = 	&\tilde{U}^0 -\frac{CFL}{360}\left(\tilde{U}^0\exp(-i4kh)-8\tilde{U}^0\exp(-i3kh)+30\tilde{U}^0\exp(-i2kh)\right)\\\nonumber&-\frac{CFL}{360}\left(-80\tilde{U}^0\exp(-ikh)+35+24\tilde{U}^0\exp(ikh)-2\tilde{U}^0\exp(i2kh)\right), \\ \nonumber
		U^{(2)} = 	&\tilde{U}^0 -\frac{CFL}{300}\left(\tilde{U}^{(1)}\exp(-i4kh)-8\tilde{U}^{(1)}\exp(-i3kh)+30\tilde{U}^{(1)}\exp(-i2kh)\right)\\\nonumber&-\frac{CFL}{300}\left(-80\tilde{U}^{(1)}\exp(-ikh)+35+24\tilde{U}^{(1)}\exp(ikh)-2\tilde{U}^{(1)}\exp(i2kh)\right), \\ \nonumber
U^{(3)} = 	&\tilde{U}^0 -\frac{CFL}{240}\left(\tilde{U}^{(2)}\exp(-i4kh)-8\tilde{U}^{(2)}\exp(-i3kh)+30\tilde{U}^{(2)}\exp(-i2kh)\right)\\\nonumber&-\frac{CFL}{240}\left(-80\tilde{U}^{(2)}\exp(-ikh)+35+24\tilde{U}^{(2)}\exp(ikh)-2\tilde{U}^{(2)}\exp(i2kh)\right), \\ \nonumber
U^{(4)} = 	&\tilde{U}^0 -\frac{CFL}{180}\left(\tilde{U}^{(3)}\exp(-i4kh)-8\tilde{U}^{(3)}\exp(-i3kh)+30\tilde{U}^{(3)}\exp(-i2kh)\right)\\\nonumber&-\frac{CFL}{180}\left(-80\tilde{U}^{(3)}\exp(-ikh)+35+24\tilde{U}^{(3)}\exp(ikh)-2\tilde{U}^{(3)}\exp(i2kh)\right), \\ \nonumber
U^{(5)} = 	&\tilde{U}^0 -\frac{CFL}{120}\left(\tilde{U}^{(4)}\exp(-i4kh)-8\tilde{U}^{(4)}\exp(-i3kh)+30\tilde{U}^{(4)}\exp(-i2kh)\right)\\\nonumber&-\frac{CFL}{120}\left(-80\tilde{U}^{(4)}\exp(-ikh)+35+24\tilde{U}^{(4)}\exp(ikh)-2\tilde{U}^{(4)}\exp(i2kh)\right), \\ \nonumber
U^{n+1} = 	&\tilde{U}^0 -\frac{CFL}{60}\left(\tilde{U}^{(5)}\exp(-i4kh)-8\tilde{U}^{(5)}\exp(-i3kh)+30\tilde{U}^{(5)}\exp(-i2kh)\right)\\\nonumber&-\frac{CFL}{60}\left(-80\tilde{U}^{(5)}\exp(-ikh)+35+24\tilde{U}^{(5)}\exp(ikh)-2\tilde{U}^{(5)}\exp(i2kh)\right). \nonumber
\end{aligned} 
\end{eqnarray}

	\begin{figure}
		\centering
	\includegraphics[width=.5\linewidth]{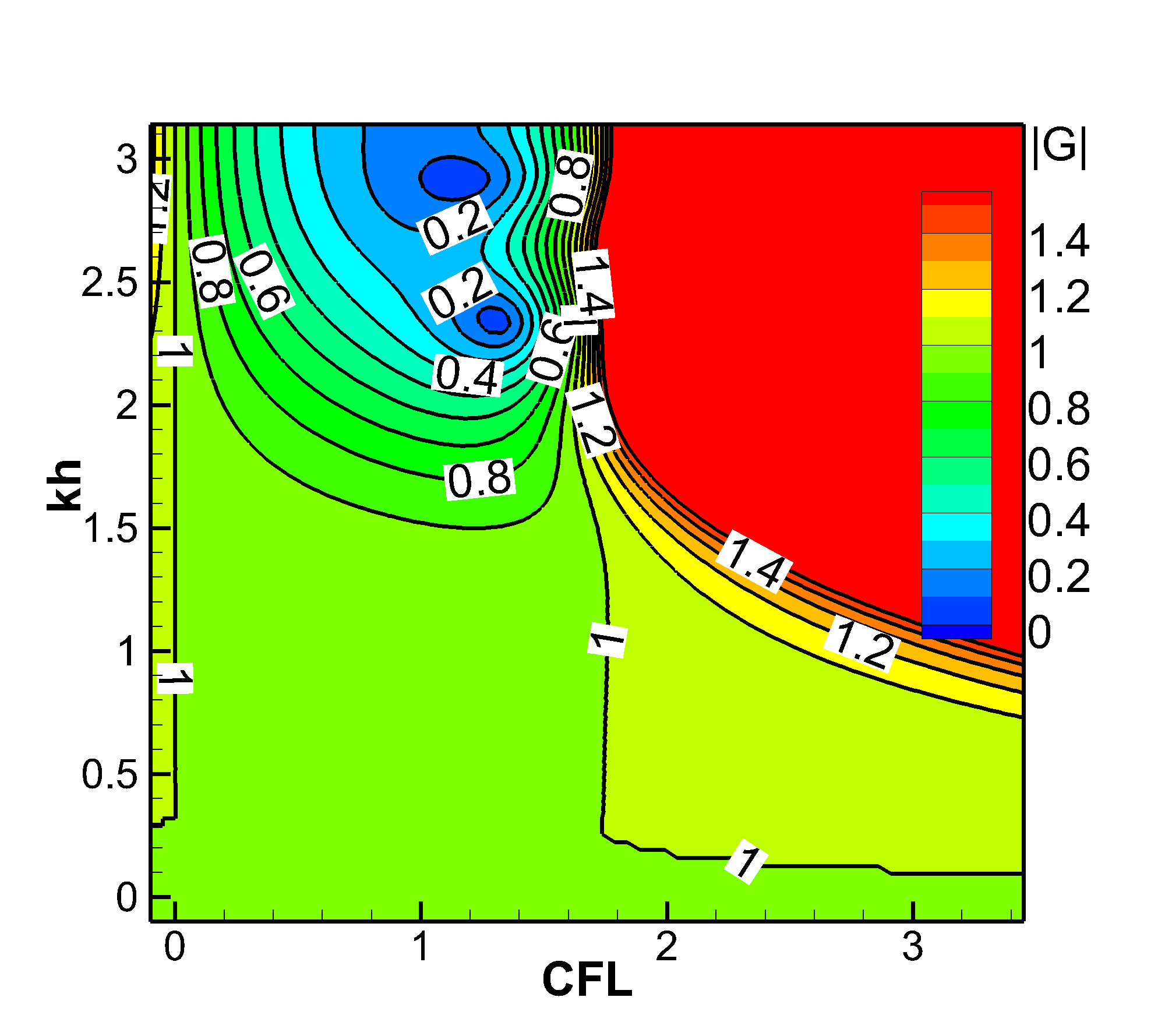}
	\caption{$|G|$ plot for RK6-L4R2 scheme }
	\label{fig:G_RK6_L4R2}
\end{figure}%

Writing in terms of numerical amplification factor gives
\begin{eqnarray}
	\begin{aligned} \label{eq:G_RK6l4r2}
		G^{(1)} = 	&1 -\frac{CFL}{360}\left(\exp(-i4kh)-8\exp(-i3kh)+30\exp(-i2kh)\right)\\\nonumber&-\frac{CFL}{360}\left(-80\exp(-ikh)+35+24\exp(ikh)-2\exp(i2kh)\right), \\ \nonumber
		G^{(2)} = 	&1 -\frac{CFL}{300}\left({G}^{(1)}\exp(-i4kh)-8{G}^{(1)}\exp(-i3kh)+30{G}^{(1)}\exp(-i2kh)\right)\\\nonumber&-\frac{CFL}{300}\left(-80{G}^{(1)}\exp(-ikh)+35+24{G}^{(1)}\exp(ikh)-2{G}^{(1)}\exp(i2kh)\right), \\ \nonumber
		G^{(3)} = 	&1 -\frac{CFL}{240}\left({G}^{(2)}\exp(-i4kh)-8{G}^{(2)}\exp(-i3kh)+30{G}^{(2)}\exp(-i2kh)\right)\\\nonumber&-\frac{CFL}{240}\left(-80{G}^{(2)}\exp(-ikh)+35+24{G}^{(2)}\exp(ikh)-2{G}^{(2)}\exp(i2kh)\right),\\ \nonumber
		G^{(4)} = 	&1 -\frac{CFL}{180}\left({G}^{(3)}\exp(-i4kh)-8{G}^{(3)}\exp(-i3kh)+30{G}^{(3)}\exp(-i2kh)\right)\\\nonumber&-\frac{CFL}{180}\left(-80{G}^{(3)}\exp(-ikh)+35+24{G}^{(3)}\exp(ikh)-2{G}^{(3)}\exp(i2kh)\right), \\ \nonumber
		G^{(5)} = 	&1 -\frac{CFL}{120}\left({G}^{(4)}\exp(-i4kh)-8{G}^{(4)}\exp(-i3kh)+30{G}^{(4)}\exp(-i2kh)\right)\\\nonumber&-\frac{CFL}{120}\left(-80{G}^{(4)}\exp(-ikh)+35+24{G}^{(4)}\exp(ikh)-2{G}^{(4)}\exp(i2kh)\right), \\ \nonumber
		G^{n+1} = 	&1 -\frac{CFL}{60}\left({G}^{(5)}\exp(-i4kh)-8{G}^{(5)}\exp(-i3kh)+30{G}^{(5)}\exp(-i2kh)\right)\\\nonumber&-\frac{CFL}{60}\left(-80{G}^{(5)}\exp(-ikh)+35+24{G}^{(5)}\exp(ikh)-2{G}^{(5)}\exp(i2kh)\right).\nonumber
	\end{aligned} 
\end{eqnarray}
Equation~\eqref{eq:G_RK6l4r2} provides the numerical gain associated with convection discretization using the RK6L4R2 method. Both theoretical analysis and simulation results yield a numerical gain of 0.9995 - i 0.0314.

From Figure~\ref{fig:G_RK6_L4R2}, it becomes evident that the scheme exhibits lower dissipation characteristics when CFL numbers are in proximity to 0 and 1.6. For situations where the CFL number falls outside of this range, it is advisable to employ a lower $kh$ value during the simulation. Numerical dissipation and dispersion tend to increase when higher values of $kh$ are used in the simulation.

\section{Numerical test cases}
\label{sec:result}
The performance of various schemes was assessed by conducting a comprehensive study on both the linear advection equation and the Navier-Stokes equation in one and two dimensions. The evaluation involved quantifying errors using the $L_1$ and $L_\infty$ norms.
\subsection{Gaussian pulse propagation}
To investigate numerical dissipation and dispersion phenomena, a Gaussian pulse with the following form
\begin{equation}\label{eq:wave_gauss}
	u(x,0) = \exp \left(\frac{-x^2}{3}\right),
\end{equation}
was advected, and the resulting numerical errors were analyzed. The computational domain [-15, 15] was discretized into 200 cells, and simulations were conducted up to time intervals of 120~s, 210~s, and 300~s. A CFL number of 0.1 was employed. 

The $L_1$ and $L_\infty$ norms, representing the discrepancies between analytical and simulated results, are summarized in table~\ref{tab:gauss}. Notably, the $L_\infty$ norm of the error for the FTBS scheme was found to be 0.6838, whereas for the RK6-L4R2 scheme, it was 5.9274E-6 with 300 mesh points. This highlights the superior performance of the RK6-L4R2 scheme, which exhibited lower numerical diffusion compared to the other schemes considered in the simulation. The results obtained for the RK4-CD4 scheme were also commendable, displaying minimal dissipation, consistent with their corresponding $G$ plots outlined in section~\ref{sec:st_adr}.
observation of the dissipation characteristics exhibited by the different schemes presented in this study. Moreover, table~\ref{tab:gauss} underscores that as the flow time increases, errors tend to grow, and the schemes become more diffusive. Consequently, it is advisable to avoid highly diffusive schemes when dealing with simulations involving extensive flow duration.

\begin{table}[h]
	\centering
	\caption{Convection equation using Gaussian pulse initial condition}
	\begin{tabular}{lcccccc}
		\toprule
		\multicolumn{1}{l}{Flow time} & \multicolumn{2}{c}{120} & \multicolumn{2}{c}{210} & \multicolumn{2}{c}{300} \\
		\midrule
		\multicolumn{1}{l}{Schemes} & \multicolumn{1}{l}{$L_1$} & \multicolumn{1}{l}{$L_{\infty}$} & \multicolumn{1}{l}{$L_1$} & \multicolumn{1}{l}{$L_{\infty}$} & \multicolumn{1}{l}{$L_1$} & \multicolumn{1}{l}{$L_{\infty}$} \\
		\midrule
		FTBS  & 24.9868 & 0.5338 & 31.57 & 0.6299 & 35.6812 & 0.6838 \\
		HRK4-CD4 & 0.0098 & 2.72E-04 & 0.0172 & 4.76E-04 & 0.0246 & 6.80E-04 \\
		RK6-L4R2 & 8.33E-05 & 2.37E-06 & 1.46E-04 & 4.15E-06 & 2.08E-04 & 5.93E-06 \\
		\bottomrule
	\end{tabular}%
	\label{tab:gauss}%
\end{table}%

\begin{figure}[ht!]
	\centering
	\includegraphics[width=0.5\linewidth]{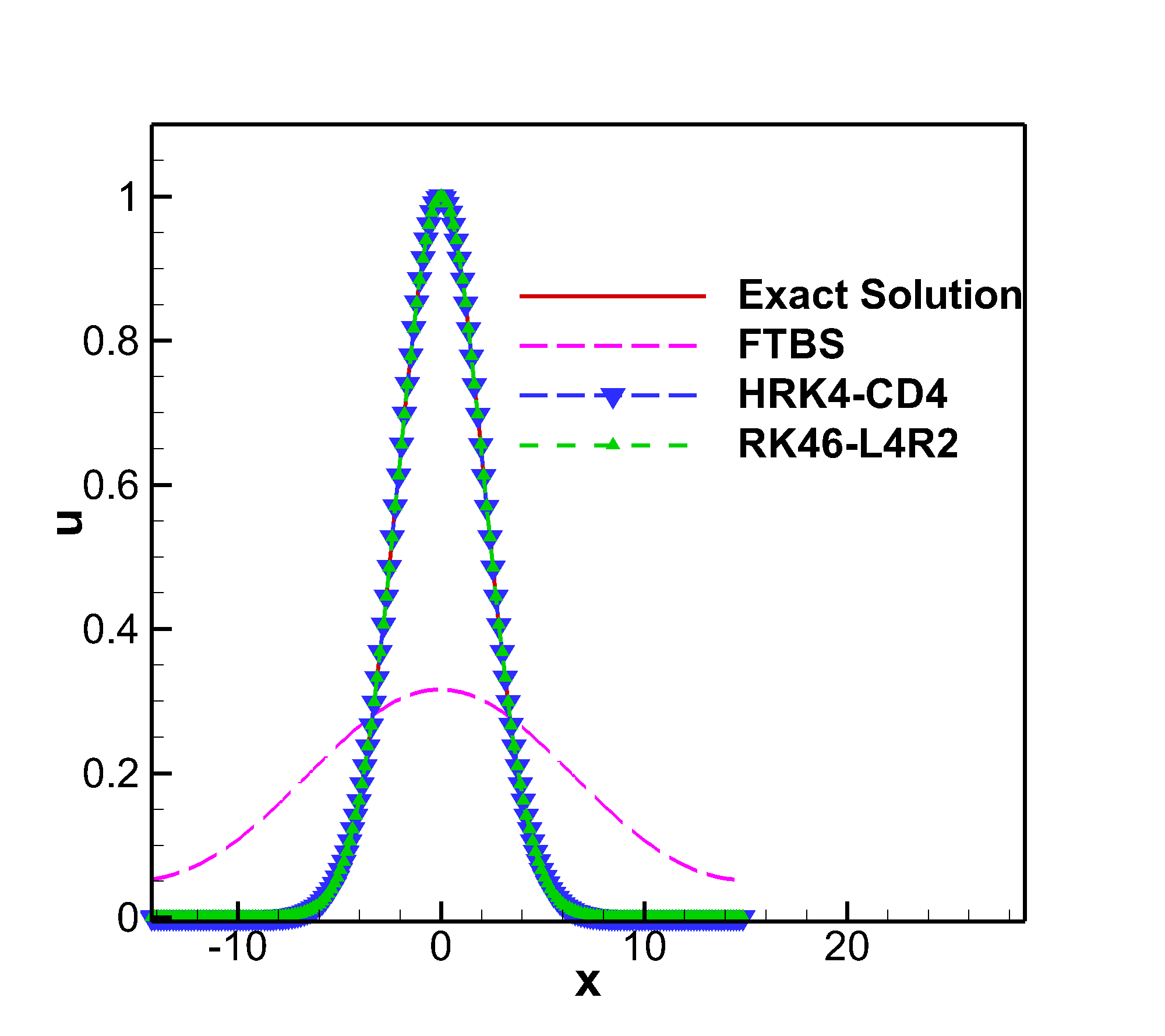}
	\caption{Solution of linear convection equation at flow time 300~s.}
	\label{fig:gaussian}
\end{figure}

\subsection{Complex wave profile simulation}
\begin{figure}[ht!]
	\begin{subfigure}[b]{0.5\textwidth}
		\includegraphics[width=\linewidth]{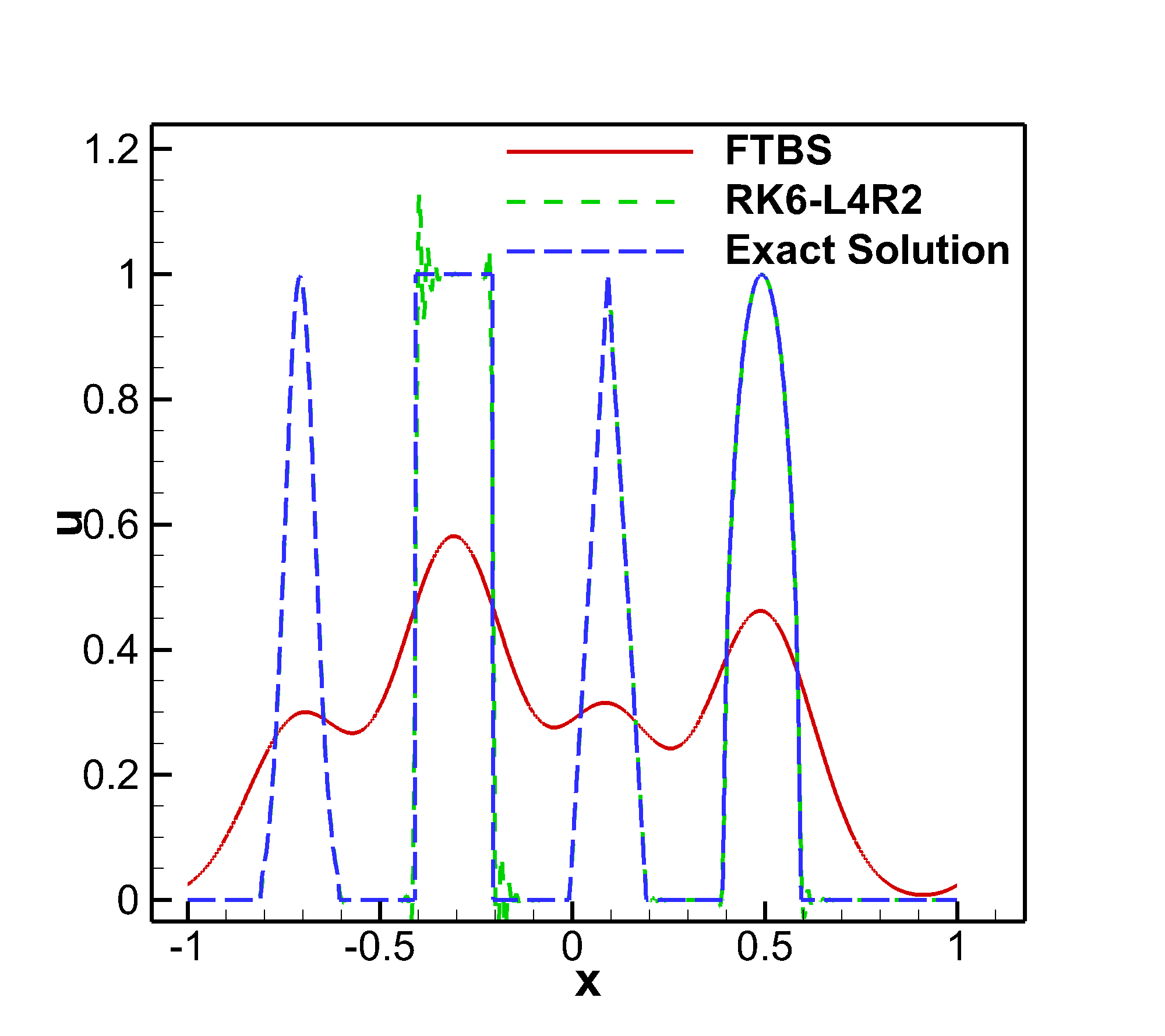}
		\caption{Solution of convection equation }
		\label{fig:Complex_ad1}
	\end{subfigure}%
	\begin{subfigure}[b]{0.5\textwidth}
		\includegraphics[width=\linewidth]{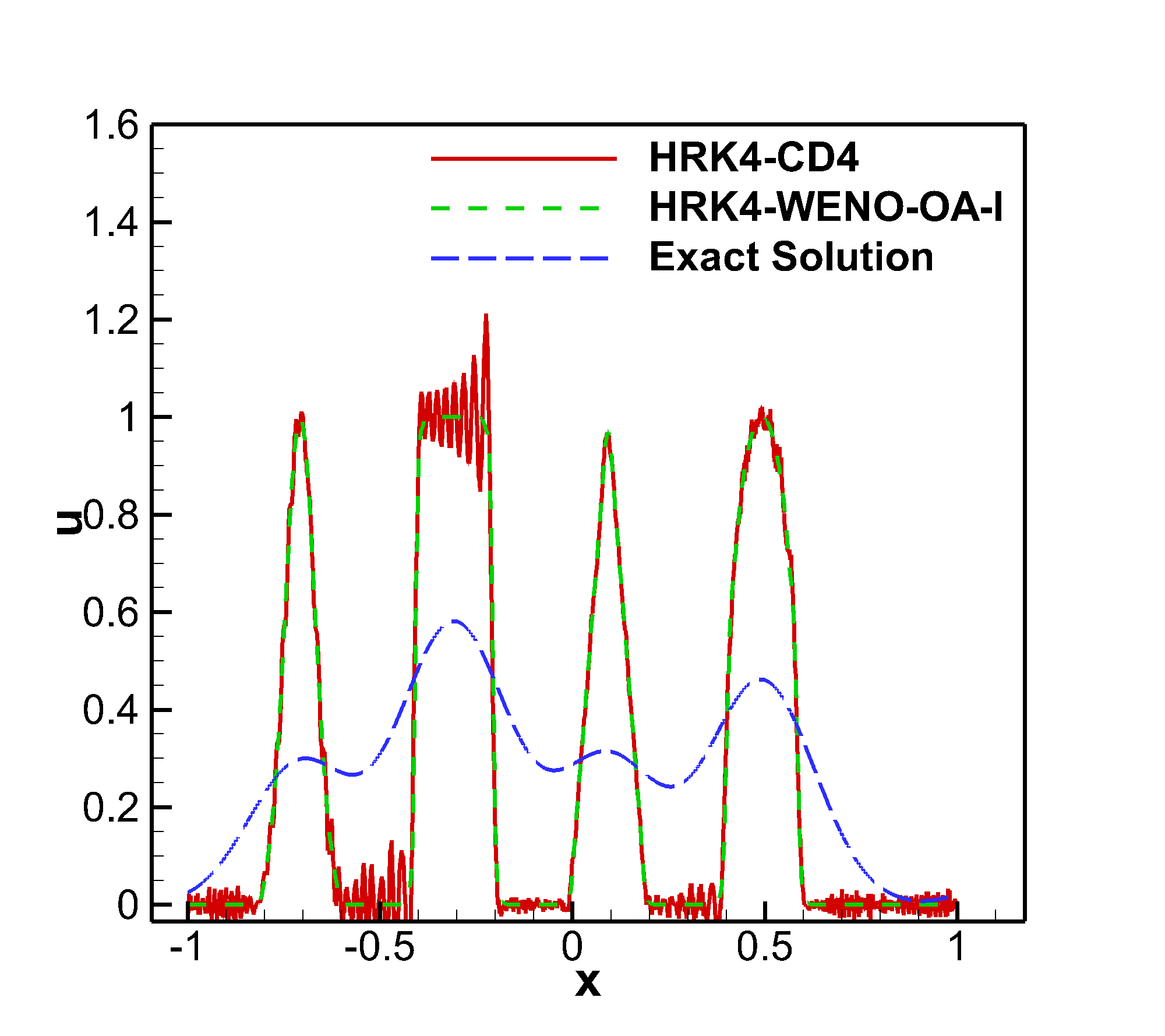}
		\caption{Solution of convection equation }
		\label{fig:Complex_ad2}
	\end{subfigure}
	\caption{Solution of convection equation at flow time 8~s using 1001 grid points using CFL = 0.25.}
	\label{fig:comp_ad}
\end{figure}

To assess the scheme's ability to handle discontinuities, the convection equation was simulated with the following initial condition~\cite{NEELAN2021100217}:
\begin{align*} 
	& u_{0}(x) =
	\begin{cases}  \displaystyle
		\frac{1}{6} \bigl(H(x, \beta, z-\delta)+H(x, \beta, z+\delta)+4 H(x, \beta, z)\bigr)  & 
		\text{if $-0.8 \leq x \leq-0.6$,}  \\[2mm]
		1 & \text{if $-0.4 \leq x \leq-0.2$,}  \\
		1-  \vert 10(x-0.1) \vert   & \text{if $0 \leq x \leq 0.2$,}  \\[2mm]
		\displaystyle
		\frac{1}{6} \bigl(L(x, \alpha, a-\delta)+L(x, \alpha, a+\delta)+4 L(x, \alpha, a)\bigr)  & 
		\text{if $0.4 \leq x \leq 0.6$,}  \\[2mm]
		0 & \text{otherwise,}
	\end{cases}
\end{align*} 
where we define the shape functions 
\begin{align*} 
	H(x, \beta, z):= \exp \bigl( -\beta(x-z)^{2} \bigr), \quad 
	L(x, \alpha, a):=\bigl( \max \{ 1-\alpha^{2}(x-a)^{2}, 0 \}\bigr)^{1/2}. 
\end{align*}
The constants are taken as $a=0.5$, $z=-0.7$, $\delta=0.005$, $\alpha=10$, and $\beta=\log 2 / (36 \delta^{2})$. The control volume $[-1, 1]$ is divided into 1001 grid points. The numerical solution is determined  up to  $T= 8~s$. Table~\ref{tab:wave_com} shows the $L_\infty$ and $L_1$ errors of the various discretization.  Figure~\ref{fig:Complex_ad1} shows the solution obtained using FTBS and RK6-L4R2 scheme. As expected FTBS scheme is diffusive. For this complex initial condition, the RK6-L4R2 scheme is dispersive. 
\begin{table}[htbp]
	\centering
	\caption{Convection equation using complex wave profile initial condition}
	\begin{tabular}{lll}
		\hline
		Schemes & \multicolumn{1}{l}{$L_1$} & $L_\infty$\\
		\hline
		FTBS  & 240.4774 & 0.6634 \\
		HRK4-CD4 & 27.2343 & 0.5352 \\
		RK6-L4R2 & 7.7960 & 0.5055 \\
		HRK41-WENO-OA-I & 8.55225 & 0.4755 \\
		\hline
	\end{tabular}%
	\label{tab:wave_com}%
\end{table}%
This discrepancy could be attributed to the fact that this scheme exhibits optimal properties only when $kh$ exceeds 0.5, and it demonstrates poor group velocity characteristics. In Figure~\ref{fig:Complex_ad2}, the results of the HRK41-WENO-OA-I scheme are displayed, showing no oscillations but a higher $L_1$ error compared to the RK6-L4R2 scheme. When employing the WENO-OA-I~\cite{neelan2023efficient} scheme to solve this problem, no oscillations appear in the solution, albeit a minor degree of dissipation is present in the sharp corners. Despite the RK6-L4R2 scheme exhibiting oscillations in this test case, it yields the lowest $L_1$ error among the considered schemes in this simulation.
\subsection{2-D convection equation}
\begin{figure}[ht!]
	\begin{subfigure}[b]{0.5\textwidth}
		\includegraphics[width=\linewidth]{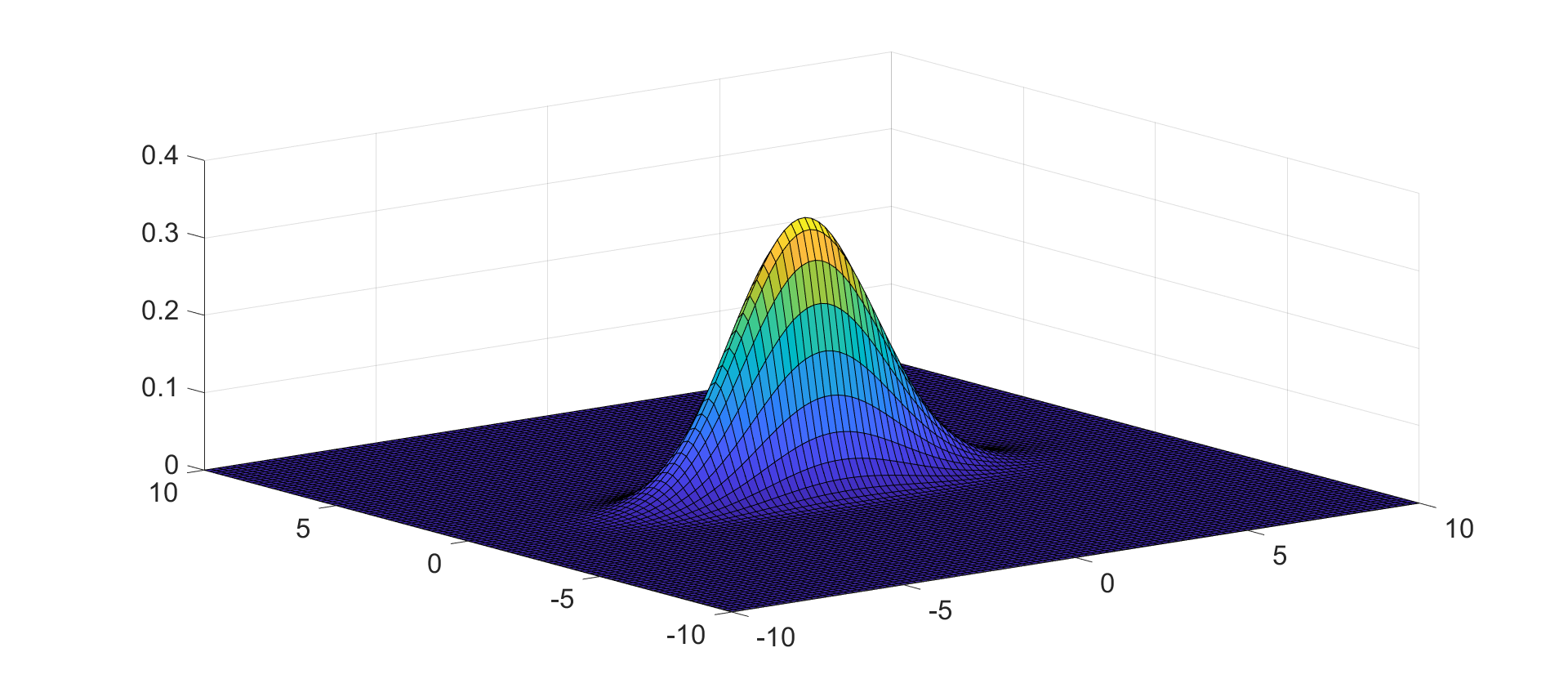}
		\caption{Using FTBS scheme }
		\label{fig:2D_FTBS}
	\end{subfigure}%
	\begin{subfigure}[b]{0.5\textwidth}
		\includegraphics[width=\linewidth]{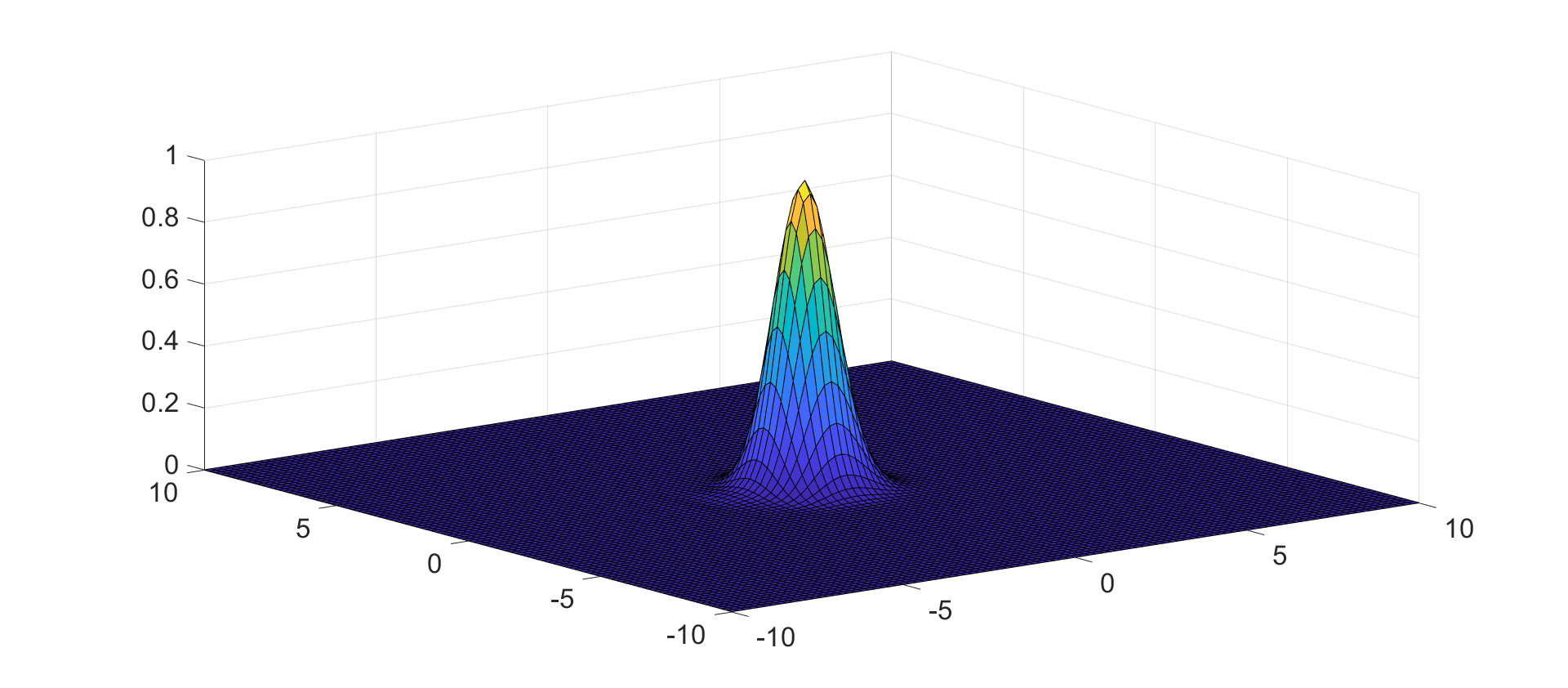}
		\caption{Using HRK4-CD4 }
		\label{fig:2D_CD4}
	\end{subfigure}
	\begin{subfigure}[b]{0.5\textwidth}
		\includegraphics[width=\linewidth]{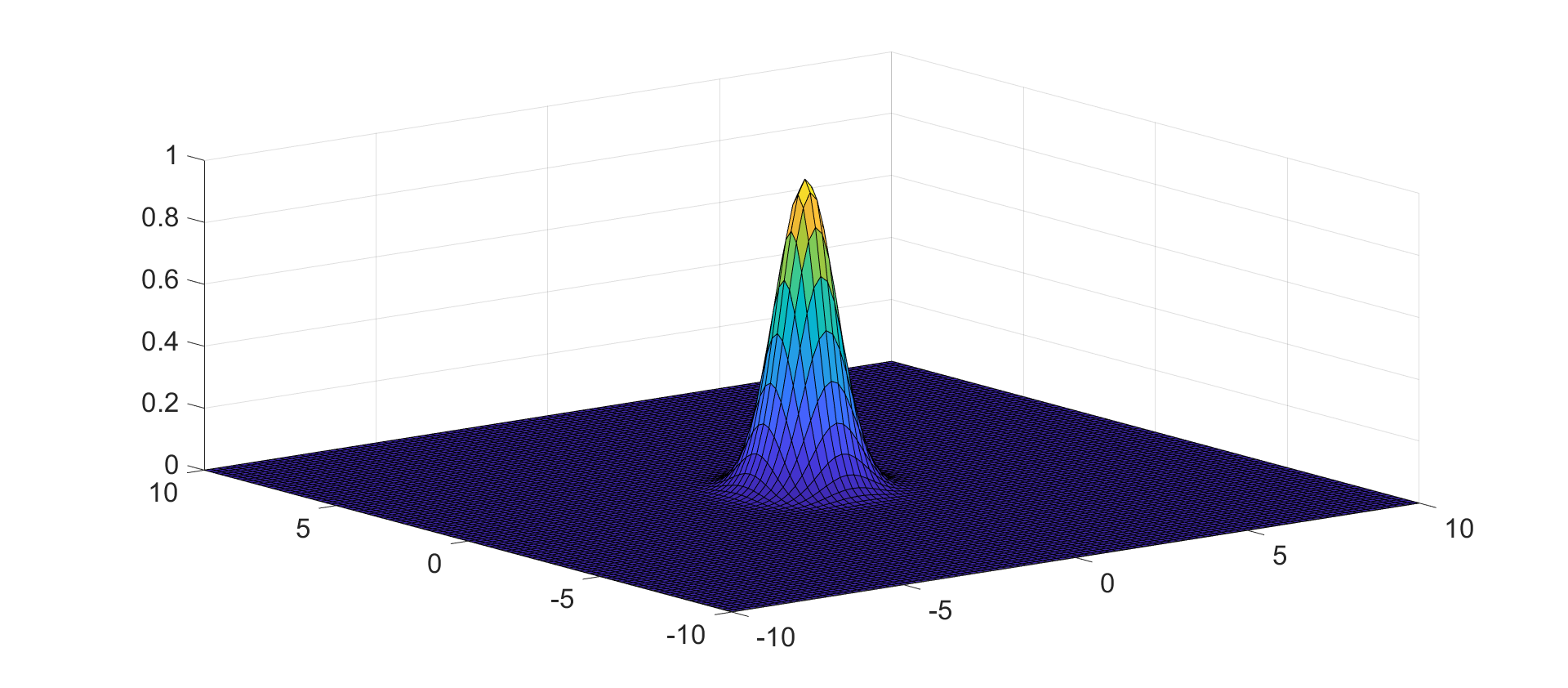}
		\caption{Using RK6-L4R2 scheme }
		\label{fig:2D_RK6L4R2}
	\end{subfigure}
	\begin{subfigure}[b]{0.4\textwidth}
		\includegraphics[width=\linewidth]{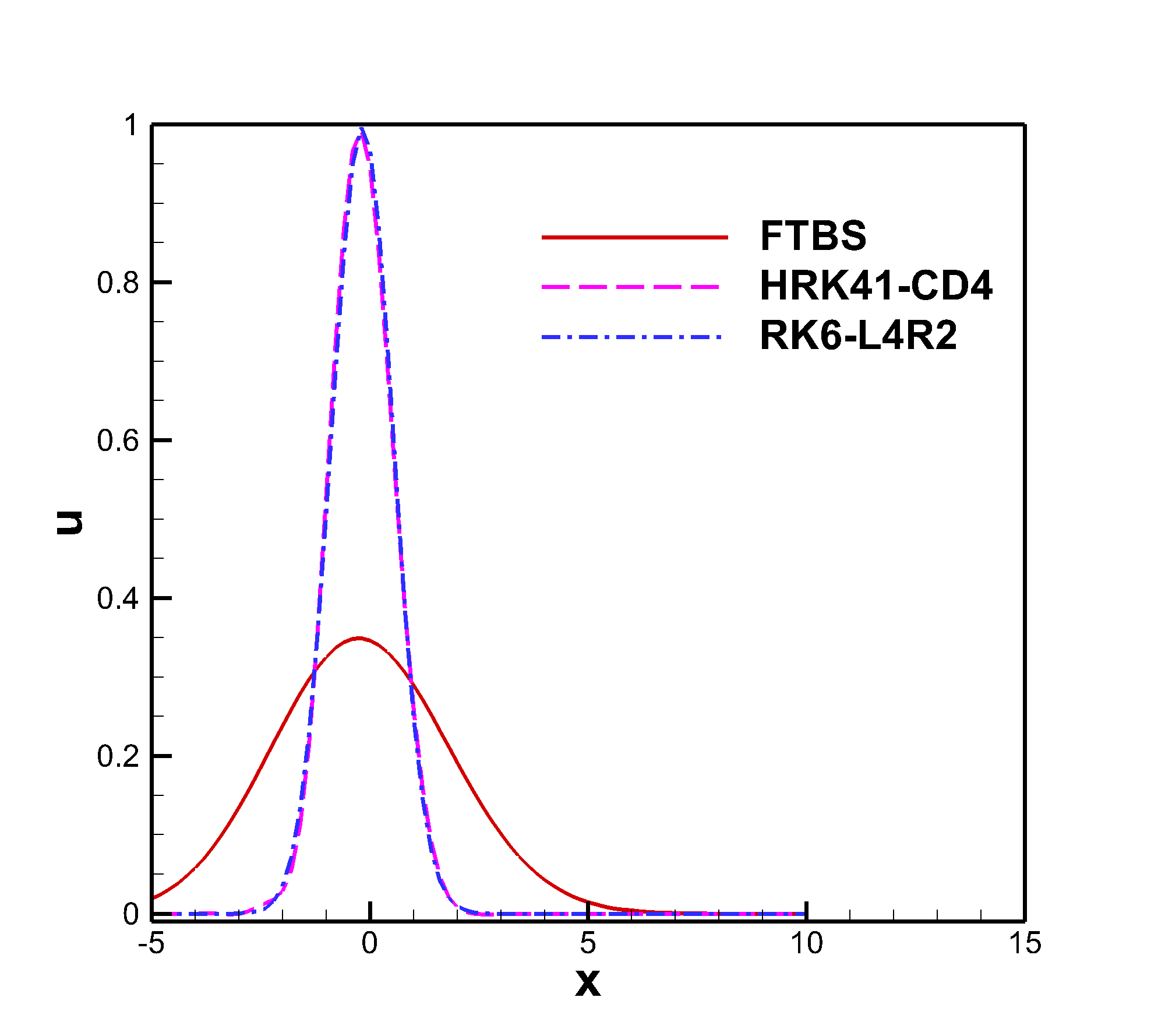}
		\caption{Solution at y = 0.5 location }
		\label{fig:Cut_sec}
	\end{subfigure}
	\caption{Solution of convection equation at flow time 20~s using 101$\times$100 grid points using CFL = 0.1.}
	\label{fig:2D_con}
\end{figure}

The governing equation is
\begin{equation}\label{key}
	\frac{\partial u}{\partial t}+c_x\frac{\partial u}{\partial x}+c_y\frac{\partial u}{\partial y}=0.
\end{equation}
The initial condition used is $u(x,0) = \exp(-x^2)\exp(-y^2)$. The domain is [-10, 10] $\times$ [-10, 10] is divided into 100 cells in each direction. We set $c_y$ to zero to easily apply periodic boundary conditions. The problem is solved up to flow time 20~s using CFL = 0.1. \begin{table}[htbp]
	\centering
	\caption{2-D Convection equation using Gaussian initial condition}
	\begin{tabular}{lll}
		\hline
		Schemes & \multicolumn{1}{l}{$L_1$} & $L_\infty$\\
		\hline
		FTBS  & 73.9927 & 0.6536 \\
		HRK4-CD4 & 16.9091 & 0.1863 \\
		RK6-L4R2 & 15.8595 & 0.1530 \\
		\hline
	\end{tabular}%
	\label{tab:2dwave}%
\end{table}
The table~\ref{tab:2dwave} shows the $L_1$ and $L_\infty$ error of the problem for different schemes. Figure~\ref{fig:2D_con} shows the solution of the 2-D convection equation using different schemes. Table~\ref{tab:2dwave} shows the error in different schemes used. Here also RK6-L4R2 scheme outperformed other schemes.

\section{Conclusion}
\label{sec:con}
A concise overview of various ADR analyses, including spatial, temporal, and space-time ADR analyses, is provided in this work. The limitations of these analyses are also discussed, and test cases are presented to illustrate these limitations. Additionally, it is noted that existing spatial-temporal ADR analyses yield the same stability equation for multi-step methods of a given order but exhibit different stability limits in numerical simulations. We have removed this limitation through the straightforward implementation of spectral analysis. 
The study includes an analysis of methods like SSPRK3, HRK41, and low-storage RK6. An optimized spatial-temporal discretization, the RK6L4R2 scheme, is also introduced, which demonstrates minimal dissipation in linear advection-based test cases. \textbf{However, it does not substantially improve results for test cases based on the Navier-Stokes equatio}n (this is not presented in the report), possibly due to its optimization for linear convection. Nonetheless, it's important to acknowledge that there are constraints to the current analysis. When dealing with nonlinear problems and non-periodic boundary conditions, a Fourier-based analysis reveals varying numerical amplification factors. An alternative approach worth considering is a deep learning-based stability assessment that concentrates on error propagation rather than the numerical amplification factor based on spectral analysis. This could potentially yield a more effective stability relationship.

\section*{Acknowledgment}
The present analysis is an improved version of the global spectral analysis presented by Prof. Tapan Sengupta in~\cite{SENGUPTA20071211}. I thank him for teaching GSA method in his course.   
\section*{Supplementary material }
Some of the Matlab codes in this work is given in the
\href{https://github.com/AGN000/Von_Neumann_Stability_matlab_code.git}{link}.
From this code, you will get an idea how to cross check the simulation result with the theoretical prediction.

\bibliography{mybibfile}.
\end{document}